# Finding Semi-Analytic Solutions of Power System Differential-Algebraic Equations for Fast Transient Stability Simulation

Nan Duan and Kai Sun, *Senior Member, IEEE*

*Abstract*— This paper studies the semi-analytic solution (SAS) of a power system's differential-algebraic equation. A SAS is a closed-form function of symbolic variables including time, the initial state and the parameters on system operating conditions, and hence able to directly give trajectories on system state variables, which are accurate for at least a certain time window. A two-stage SAS-based approach for fast transient stability simulation is proposed, which offline derives the SAS by the Adomian Decomposition Method and online evaluates the SAS for each of sequential time windows until making up a desired simulation period. When applied to fault simulation, the new approach employs numerical integration only for the fault-on period to determine the post-disturbance initial state of the SAS. The paper further analyzes the maximum length of a time window for a SAS to keep its accuracy, and accordingly, introduces a divergence indicator for adaptive time windows. The proposed SAS-based new approach is validated on the IEEE 10-machine, 39-bus system.

*Index Terms*—Differential-algebraic equations; High performance computing; Laplace Adomian Decomposition; semi-analytic solution; transient stability; time-domain simulation

## I. Introduction

TIME-domain simulation of a power system following a contingency for transient stability analysis needs to solve nonlinear differential-algebraic equations (DAEs) on the system state over a simulation period. Numerical integration methods, either explicit or implicit, are usually employed to solve the Initial Value Problem (IVP) of the DAEs but could be time-consuming for a large-scale power system with many generators because the DAEs in nature model tight coupling between generators via nonlinear sine functions. Also, numerical instability is a concern with explicit integration methods like the Runge–Kutta method, which is widely applied in todays' simulation software; implicit integration methods like the Trapezoidal method overcome numerical instability by introducing implicit algebraic equations, which also need to be solved by numerical methods like the Newton-Raphson method, and thus, the computational complexity is significantly increased.

Intuitively, if the analytical solution of the IVP of power system DAEs could be found as a set of explicit algebraic functions of symbolic variables including time, the initial state and parameters on the system operating condition, it would directly give the values of state variables at any time without conducting time-consuming iterations as numerical integration does. However, such an analytical solution which is accurate for any simulation time period does not exist in theory for power system DAEs. Thus, a compromise could be to find an approximate analytic solution, named a semi-analytic solution (SAS), which keeps accuracy for a certain time window, denoted by $T$, and repeats using the SAS until making up a desired simulation period. Thus, solving the IVP becomes simply evaluating the SAS, i.e. plugging in values of symbolic variables, and can be extremely fast compared to numerical integration. If evaluation of the SAS over each window $T$ takes a time less than $\tau$, the IVP is solved will be $T/\tau$ times faster than real time.

The true solution of power system DAEs may be approached by summating infinite terms of some series expansion. A SAS can be defined as the sum of finite terms that is accurate over window $T$. Such a series expansion can be derived by using the Adomian Decomposition Method (ADM). Compare to other decomposition methods like Taylor's Series Expansion, The ADM is able to keep nonlinearity in the system model [1]. In this paper, the Multi-stage Modified ADM (MM-ADM) is applied to derive a SAS of the classical 2$^{nd}$ order power system DAE model for transient stability simulation. The rest of the paper is organized as follows. Section II introduces the fundamental concepts of the ADM to MM-ADM. Then section III illustrates how a SAS is derived for a single-machine-infinite-bus (SMIB) system. The maximum window of accuracy of the SAS, i.e. the limit of $T$, is studied. Section IV proposes a new approach to fast transient stability simulation by offline deriving and online evaluating a SAS. In section V, the proposed approach is tested on the IEEE 10-machine, 39-bus system. Finally, conclusions are drawn in section VI.

## II. Introduction of the ADM and MM-ADM

This section briefly introduces the concepts of the ADM, the modified ADM and the multi-stage ADM [1]-[5].

### A. ADM and Modified ADM

As presented later in Section III, the classical 2$^{nd}$ order DAE model for a SMIB power system can be transformed into a 2$^{nd}$ order differential equation in the form of (1), where $f$ is a

N. Duan and K. Sun are with the department of Electrical Engineering and Computer Science, University of Tennessee, Knoxville, TN 37996 USA (e-mail: nduan@vols.utk.edu, kaisun@utk.edu).
The work is supported by the University of Tennessee, Knoxville.

nonlinear differentiable function and a is a coefficient related to the system's oscillation damping.

$$\ddot{x}(t) + a\dot{x}(t) = f(x(t)) \quad (1)$$

To solve its IVP, the first step of ADM is to apply Laplace transform $L[\cdot]$ to both sides to transform the differential equation about time $t$ into an algebraic equation about complex frequency $s$, i.e. (2). Then solve $L[x]$ to get (3).

$$s^2 L[x] - x(0)s - \dot{x}(0) + aL[\dot{x}] = L[f(x)] \quad (2)$$

$$L[x] = \frac{x(0)}{s} + \frac{\dot{x}(0)}{s^2} - \frac{aL[\dot{x}]}{s^2} + \frac{L[f(x)]}{s^2} \quad (3)$$

Assume that the solution $x(t)$ can be decomposed as follows

$$x(t) = \sum_{n=0}^{\infty} x_n(t) \quad (4)$$

Decompose $f(x)$ to be the sum of Adomian polynomials calculated from (6) by introducing a Lagrange multiplier $\lambda$

$$f(x) = \sum_{n=0}^{\infty} A_n(x_0, x_1, \cdots, x_n) \quad (5)$$

$$A_n = \frac{1}{n!}\left[\frac{\partial^n}{\partial \lambda^n} f\left(\sum_{i=0}^{n} x_i \lambda^i\right)\right]_{\lambda=0} \quad (6)$$

For instance, the first five Adomian polynomials are following, where derivatives are with respect to $\lambda$ and take values at $\lambda=0$:

$$A_0 = f(x_0) \quad (7)$$

$$A_1 = x_1 f'(x_0) \quad (8)$$

$$A_2 = x_2 f'(x_0) + \frac{x_1^2}{2!} f''(x_0) \quad (9)$$

$$A_3 = x_3 f'(x_0) + x_1 x_2 f''(x_0) + \frac{x_1^3}{3!} f'''(x_0) \quad (10)$$

$$A_4 = x_4 f'(x_0) + \left(x_1 x_3 + \frac{x_2^2}{2}\right)f''(x_0) + \frac{x_1^2 x_2}{2} f^{(3)}(x_0) + \frac{x_1^4}{4!} f^{(4)}(x_0) \quad (11)$$

There are also alternative methods to calculate Adomian polynomials. A more computationally efficient algorithm for calculating $A_0, A_1, \ldots, A_n$ introduced by [2] is applied in this paper as follows:
*Step-1*: Set $A_0 = f(x_0)$.
*Step-2*: For $k=0$ to $n-1$, conduct

$$A_k(x_0, \ldots, x_k) := A_k(u_0 + u_1\lambda, \ldots, u_k + (k+1)u_{k+1}\lambda)$$

*Step-3*: Calculate $A_{k+1} = \frac{1}{k+1}\frac{\partial}{\partial \lambda} A_k \Big|_{\lambda=0}$.

After decomposing both $x(t)$ and $f(x)$, compare their terms to easily derive these recursive formulas for $L[x_n]$.

$$L[x_0] = \frac{x(0)}{s} + \frac{\dot{x}(0)}{s^2} \quad (12)$$

$$L[x_{n+1}] = \frac{aL[\dot{x}_n]}{s^2} + \frac{L[A_n]}{s^2}, n \geq 0 \quad (13)$$

By applying inverse Laplace transform to both sides of (12) and (13), $x_n$ for any $n$ can be obtained. Consequently, a SAS of (1) is yielded by summating first $N$ terms of $x_n$:

$$x^{SAS}(t) = \sum_{n=0}^{N-1} x_n(t) \quad (14)$$

Computational burden is caused when $x_n$ is plugged into (6) recursively to calculate $A_n$. Since in formula (6), $x_0$ always appears in $f$ and its derivatives, this paper adopts a modified ADM proposed by [14] to calculate the following (15)-(17) instead of (12) and (13) to reduce the computational burden. Thus, the SAS becomes a power series expansion with only polynomial nonlinearity. However, the calculation of the approximate solution is significantly accelerated.

$$L[x_0] = \frac{x(0)}{s} \quad (15)$$

$$L[x_1] = \frac{\dot{x}(0)}{s^2} + \frac{aL[\dot{x}_0]}{s^2} + \frac{L[A_0]}{s^2} \quad (16)$$

$$L[x_{n+2}] = \frac{aL[\dot{x}_{n+1}]}{s^2} + \frac{L[A_{n+1}]}{s^2}, n \geq 0 \quad (17)$$

B. *For a multi-variable system*

Similarly, the classical 2nd order DAE model of a $K$-machine power system can be transformed into a $K$-variable system whose each variable follows a 2nd order differential equation coupled with the other variables through a nonlinear function. To solve the $K$ differential equations, the ADM works in the same manner as for the single-variable system except for computation of Adomian polynomials regarding each variable. Consider a $K$-variable system given by (18) with $K$ nonlinear functions $f_1, f_2 \ldots f_K$ which are respectively decomposed by (20).

$$\begin{aligned}
\ddot{\mathbf{x}}(t) + \mathbf{a}\dot{\mathbf{x}}(t) &= \mathbf{f}(\mathbf{x}(t)) \\
\mathbf{x}(t) &= \begin{bmatrix} x_1(t) & x_2(t) & \cdots & x_K(t) \end{bmatrix}^T \\
\mathbf{f}(\cdot) &= \begin{bmatrix} f_1(\cdot) & f_2(\cdot) & \cdots & f_K(\cdot) \end{bmatrix}^T \\
\mathbf{a} &= \begin{bmatrix} a_1 & 0 & 0 & 0 \\ 0 & a_2 & 0 & 0 \\ 0 & 0 & \ddots & 0 \\ 0 & 0 & 0 & a_K \end{bmatrix}
\end{aligned} \quad (18)$$

$$\mathbf{x}(t) = \sum_{n=0}^{\infty} \mathbf{x}_n(t) \quad (19)$$

$$\mathbf{x}_n = \begin{bmatrix} x_{1,n} & x_{2,n} & \cdots & x_{K,n} \end{bmatrix}^T$$

$$f_i(\mathbf{x}) = \sum_{n=0}^{\infty} A_{i,n}(\mathbf{x}_0, \mathbf{x}_1, \ldots, \mathbf{x}_n), \quad i = 1 \cdots K \quad (20)$$

where Adomian polynomials $A_{i,n}$'s are respectively calculated by (21). Then, SAS's of $x_1, x_2 \ldots x_K$ can be obtained similarly.

$$A_{i,n} = \frac{1}{n!}\left[\frac{\partial^n}{\partial \lambda^n} f_i\left(\sum_{i=0}^{n} \mathbf{x}_i \lambda^i\right)\right]_{\lambda=0} \quad (21)$$

C. *Multi-stage ADM*

Since the SAS of the IVP is accurate only for a limited time window $T$, a Multi-stage ADM is adopted to extend accuracy to an expected simulation period [7][10][13] by two steps:

*Step-1:* Partition the simulation period into sequential windows of $T$ which each can keep an acceptable accuracy of the SAS.

*Step-2:* Evaluate the SAS for the first $T$ using a given initial state and the values of other parameters; starting from the second $T$, evaluate the SAS by taking the final state of the previous $T$ as the initial state.

This paper utilizes the above Multi-stage Modified ADM (MM-ADM) for transient stability simulation.

## III. APPLICATION TO A SMIB SYSTEM

### A. Deriving a SAS by the Modified ADM

This paper focuses on using the SAS-based approach for transient stability analysis over several seconds following a disturbance. Thus, the 2nd order classical generator model can meet that requirement and is considered in the rest of the paper. To illustrate the details on how a SAS is derived, a SMIB system modeled by DAE (22) is studied in this section:

$$\begin{cases} \dfrac{2H}{\omega_0} \dfrac{d\Delta\omega}{dt} + D\Delta\omega = P_m - P_e \\ \dfrac{d\Delta\delta}{dt} = \Delta\omega \\ P_e = P_e(\Delta\delta) = E^2 G + EE_\infty Y \cos(\Delta\delta + \delta_0 - \theta - \delta_\infty) \end{cases} \quad (22)$$

$\omega_0$ and $\Delta\omega$ are the synchronous speed and the deviation of the rotor speed. $\delta_0$ and $\Delta\delta$ are respectively the steady-state value and deviation of the rotor angle. $H$ and $D$ are respectively the inertia and damping coefficient of the generator. $E$ is EMF magnitudes of the generator and the infinite bus. $Y\angle\theta$ is the admittance between the generator and the infinite bus. $G$ is the self-conductance of the generator. $E_\infty\angle\delta_\infty$ is the voltage phasor of the infinite bus, and usually let $\delta_\infty$ be 0. $P_e$ is the instantaneous electrical power output changing with $\Delta\delta$ and $P_m=P_e(0)$ is the mechanical power input representing the operating condition, i.e. the steady-state value of $P_e$. Rewrite (22) in the form of (1):

$$\Delta\ddot{\delta} + \dfrac{D}{2H}\Delta\dot{\delta} = \dfrac{\omega_0}{2H}\left[ P_m - E^2 G - EE_\infty Y \cos(\delta_0 + \Delta\delta - \theta) \right] \quad (23)$$

Assume

$$\Delta\delta(t) \approx \sum_{n=0}^{N-1} \Delta\delta_n(t) \quad (24)$$

Apply (6) to derive $A_n$. For example, first 3 terms are

$$A_0 = \dfrac{\omega_0}{2H}\left[ P_m - E^2 G - EE_\infty Y \cos(\delta_0 + \Delta\delta_0 - \theta) \right] \quad (25)$$

$$A_1 = -\dfrac{\omega_0 EE_\infty Y}{2H} \Delta\delta_1 \sin(\delta_0 + \Delta\delta_0 - \theta) \quad (26)$$

$$A_2 = \dfrac{\omega_0 EE_\infty Y}{2H} [\dfrac{1}{2}\Delta\delta_1^2 \cos(\delta_0 + \Delta\delta_0 - \theta) + \Delta\delta_2 \sin(\delta_0 + \Delta\delta_0 - \theta)] \quad (27)$$

TABLE I
Parameters of the SMIB System

| $D$ | $H$ | $EE_\infty Y$ | $G$ |
|---|---|---|---|
| 1pu | 3s | 1.7pu | 0pu |

| $\delta_0$ | $\omega_0$ | $\Delta\delta(0)$ | $\Delta\dot{\delta}(0)$ |
|---|---|---|---|
| 1.0472rad | 377rad/s | 0.0957rad | 3.7639rad |

$D$, $H$, $Y$, $E$, $E_\infty$ and $G$ are parameters about the system model, $\delta_0$ and $\omega_0$ determine the operating condition, and $\Delta\delta(0)$ and $\Delta\dot{\delta}(0) = \Delta\omega(0)$ give the post-disturbance initial state. Any of these parameters can be treated as symbolic variables together with time in the SAS. To illustrate the derivation of the SAS, consider a simple case that only symbolizes time and has the other parameters given in Table I. Follow the procedure of (15)-(17) to derive first 5 terms:

$$\Delta\delta_0 = 0.0957 \quad (28)$$

$$\Delta\delta_1 = -2.3400t^2 + 3.7639t \quad (29)$$

$$\Delta\delta_2 = 8.6433t^4 - 27.6760t^3 - 0.3137t^2 \quad (30)$$

$$\Delta\delta_3 = -3.9015t^6 + 18.2507t^5 + 59.6802t^4 + 0.01743t^3 \quad (31)$$

$$\Delta\delta_4 = -33.7017t^8 + 216.8390t^7 - 448.2592t^6 + 11.9317t^5 - 0.0007t^4 \quad (32)$$

Summating them, a SAS with $N=5$ of $\Delta\delta(t)$ is:

$$\Delta\delta(t) = -33.7017t^8 + 216.8390t^7 - 452.1607t^6 + 30.1824t^5 \\ + 68.3227t^4 - 27.6585t^3 - 2.6536t^2 + 3.7639t + 0.0957 \quad (33)$$

Similarly, the SAS's for any $N$ can be derived. Fig. 1 compares the curves of SAS's with $N=5\sim8$ with the true solution estimated by the 4th order Runge-Kutta (RK4). Those SAS curves match well the true solution for first 0.2-0.4s depending on $N$. The bigger $N$, the longer period of accuracy.

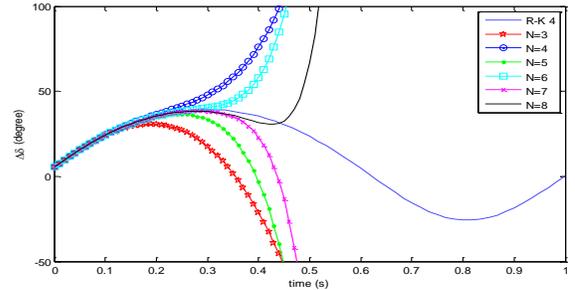

Fig. 1. Comparison of SASs with numerical result

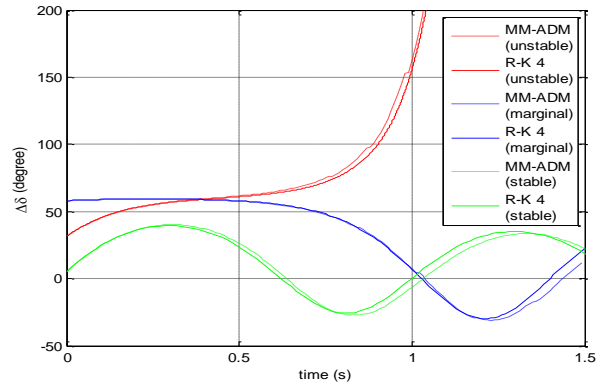

Fig. 2. Comparisons for three cases

Adjust $\Delta\delta(0)$ and $\Delta\dot{\delta}(0)$ to create typical stable, marginally stable and unstable cases. The MM-ADM is tested to derive a SAS with $N=3$ and repeat its evaluation every 0.17s. The results compared to RK4 are shown in Fig. 2. All the three SAS curves match true solutions well.

## B. Maximum Window of Accuracy

To study the maximum window keeping the accuracy of a SAS, i.e. the limit of $T$, an indicator for a SAS to lose its accuracy is proposed and illustrated using the SMIB system. For example, consider the SAS with $N=5$ given in (33). The curves of the SAS and its first 5 terms are drawn in Fig. 3 to compare with the true solution estimated by RK4.

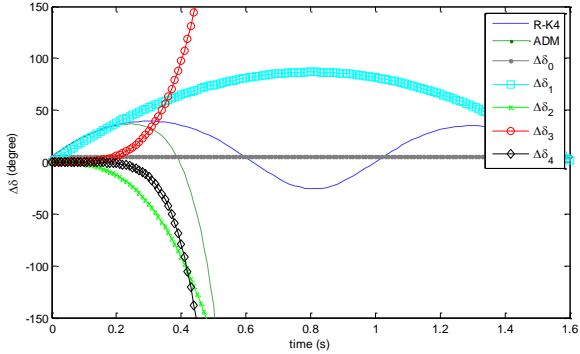

Fig. 3. Comparison between each decomposition terms of N=5 result

There are two observations from the figure. First, the SAS matches well the true solution within a time window (about 0.2s for this case), which is defined as the maximum window of accuracy and denoted by $R_A$, beyond which $\Delta\delta_n$ terms with a large $n$ lead to infinite quickly. Second, the larger $n$, the smaller contribution of term $\Delta\delta_n$ to accuracy of the SAS. Accordingly, the derivative of the last $n$ term, i.e. $d\Delta\delta_{N+1}/dt$, is defined as an indicator of loss of accuracy, denoted by $I_{LOA}$, which is close to zero within $R_A$ and sharply increases, otherwise. By selecting an appropriate threshold $I_{LOA,max}$ for $I_{LOA}$, $R_A$ can be calculated by means of the formula of a SAS. In the following, how $R_A$ depends on time constants of the SMIB system is analyzed and the conclusion drawn is then tested on the IEEE 3-machine, 9-bus system.

Given a SAS with $N$ terms for the SMIB system, the relationship between $R_A$ and the inertia $H$ can be developed with $I_{LOA,max}$ defined. To make this formula more general, relax the value of $\delta_\infty$ in (22), which in a multi-machine system could represent the instantaneous rotor angle of a remote large machine. If $N=3$, the 3rd term of the SAS for generator is

$$\delta_3 = c_1 t^4 + c_2 t^3$$

$$c_1 = \frac{\omega_0^2 Y E E_\infty \sin(\theta + \delta_\infty(0) - \delta(0))}{96 H^2} \{(E^2 G - P_m) \quad (34)$$
$$+ Y E E_\infty [\cos(\theta + \delta_\infty(0) - \delta(0))]\}$$

$$c_2 = \frac{\omega_0 Y E E_\infty [\dot\delta_\infty(0) - \dot\delta(0)] \sin(\theta + \delta_\infty(0) - \delta(0))}{12 H}$$

$\delta(0)$ and $\dot\delta(0)$ are the post-fault initial value and initial derivative of generator's power angle. Take the first order derivative of (34) with respect to $t$, then substitute $t$ with $R_A$ in and let $\dot\delta_3 = I_{LOA,max}$ to derive

$$I_{LOA,max} = 4c_1 R_A^3 + 3c_2 R_A^2 \quad (35)$$

Although (35) is derived from the SMIB system, it can also be applied to a multi-machine system to approximately estimate $R_A$ for each machine. For that case, $E_\infty \angle \delta_\infty$ could be the voltage phasor of a reference node with basically consistent voltage magnitude and angle, e.g. the slack bus in the power-flow model, an equivalent bus of a neighboring system and the EMF of a large machine. Estimate the transfer admittance $Y \angle \theta$ from that machine to the reference node by eliminating the other rows and columns in the system admittance matrix containing all network nodes and generator internal nodes.

Formulas (35) gives how $R_A$ depends on the inertia, i.e. $H$, and the initial state, i.e. $\delta(0)$ and $\dot\delta(0)$. For any given parameters of the SMIB system, $R_A$ can be estimated as the upper limit of $T$ to perform MM-ADM. On the contrary, if a desired $R_A$ is provided, (35) can also estimate the expected minimum inertia of the machine, denoted by $H_{min}$, to achieve that $R_A$. That is helpful information for dynamic model reduction of a large-scale power system for the purpose of speeding up simulation: generators with small inertias can be aggregated by coherent groups or geographic regions to equivalent generators that have inertias larger than $H_{min}$.

It should be pointed out that because (35) ignores, or in other words simplifies the coupling between that studied machine with the other machines, if there is any other machine with a small inertia along the path from that studied machine to the reference node, it will influence the actual $R_A$ or the lower limit of $H$ of the studied machine. Thus, a safer approach could be first estimating $R_A$ (or equivalently, $H_{min}$) for each machine using (35), and then choosing the minimum $R_A$ (or the maximum $H_{min}$) for all generators.

The IEEE 3-generator 9-bus system in [6] is studied to validate (35) for a multi-machine system. As shown in Table II, gradually decrease $H_3$, the inertia of generator 3, from 4.5s to 1.0s and keep the other two unchanged at original 23.64s and 6.4s, such that 8 system models are yielded. For each model, $R_A$ is estimated and given in the table. Taking $H_3=3.0s$ as an example, bus 1 is selected as the reference node and has voltage equal to 1.0170 pu. The transfer admittance from buses 1 to the EMF of generator 3 is $Y = 1.0792 \angle 80.27° pu$. Choose $I_{LOA,max}$ as 5 rad/s. Since two oscillation modes exist in this system, their oscillation periods $T_1$ and $T_2$ are important time constants influencing $R_A$, which are calculated by eigen-analysis on each linearized system model as listed in Table II.

TABLE II
$R_A$ VS. TIME CONSTANTS OF THE SYSTEM

| No. | $H_3$(s) | $R_A$(s) | $T_1$(s) | $T_2$(s) |
|---|---|---|---|---|
| 1 | 4.5 | 0.2546 | 0.9510 | 0.5516 |
| 2 | 4.0 | 0.2342 | 0.9438 | 0.5280 |
| 3 | 3.5 | 0.2131 | 0.9369 | 0.5014 |
| 4 | 3.0 | 0.1905 | 0.9304 | 0.4718 |
| 5 | 2.5 | 0.1662 | 0.9241 | 0.4365 |
| 6 | 2.0 | 0.1410 | 0.9183 | 0.3961 |
| 7 | 1.5 | 0.1137 | 0.9128 | 0.3479 |
| 8 | 1.0 | 0.0845 | 0.9076 | 0.2881 |

The same contingency which cuts the line between bus 5 and bus 7 during fault and closes the line after fault is applied to each case with different $H_3$ parameters.

From the table, the bigger time constant $T_1$ does not change significantly with $H_3$, so Fig. 4 only draws how $R_A$ changes with the smaller time constant $T_2$. The relationship is monotonic and basically linear. A conclusion to draw is that for a system with multiple oscillation modes, $R_A$ is mainly influenced by the fastest mode, which relates to generators with small inertias. Thus, for a multi-machine system with small generators, if instability with small generators is not of interests in transient stability analysis, those small generators can be merged with their neighboring generators to eliminate the inertias smaller than the lower limit of $H$ for achieving a given $R_A$.

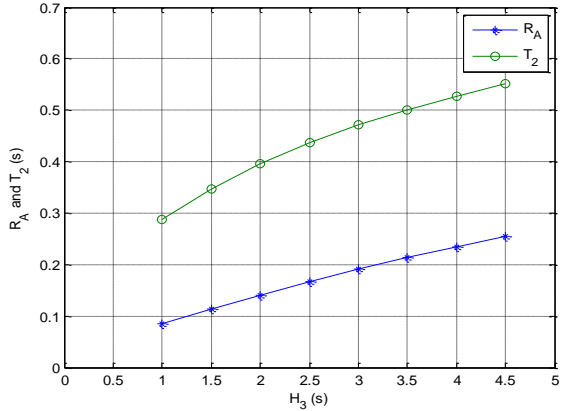

Fig. 4. Relationships between $R_A$, $T_2$ and $H_3$

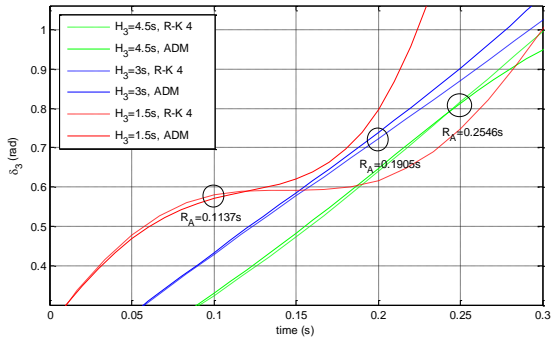

Fig. 5. $R_A$'s with respect to selected $H_3$'s

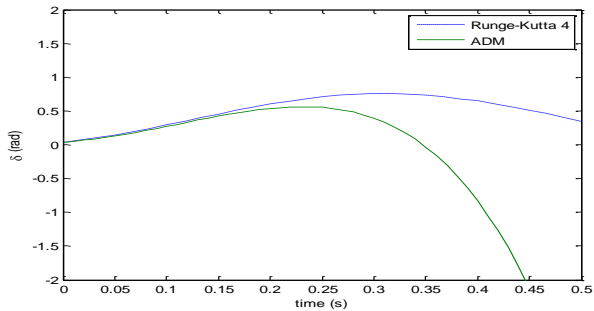

Fig. 6. Maximum $\dot{\delta}(0)$ starting point

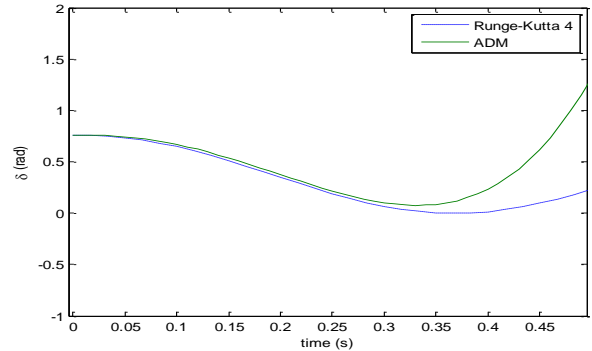

Fig. 7. Minimum $\dot{\delta}(0)$ starting point

In Fig. 5, the comparison results of Modified ADM and R-K 4 are given for 3 selected conditions with different $H_3$'s. As shown in Fig. 5, these $R_A$'s agree with the estimations obtained from (35). The dependency of $R_A$ on the initial state $\delta(0)$ and $\dot{\delta}(0)$ is illustrated in Fig. 6 and Fig. 7: $R_A$ is around 0.1s if the evaluation of the SAS starts from an initial state with large $|\dot{\delta}(0)|$; $R_A$ increases to around 0.25s when $|\dot{\delta}(0)|$ decreases to 0. A conclusion is that $R_A$ may be extended by appropriately choosing the starting point when evaluating the same SAS. Thus, the SAS-based approach may cooperate with the numerical integration approach: the numerical approach gives the trajectory from the beginning of the disturbance to the first time with $|\dot{\delta}(0)|=0$ and then the SAS will provide the trajectory thereafter. A SAS's dependency on the initial state also indicates that $R_A$ should be estimated for each time of SAS evaluation to allow an adaptive $T$ during the simulation period for the best accuracy. Alternatively, a fixed $T$ may also be used to save the time of $R_A$ evaluation as long as it is smaller than $R_A$ with any credible initial state.

## IV. PROPOSED MM-ADM BASED TWO-STAGE APPROACH

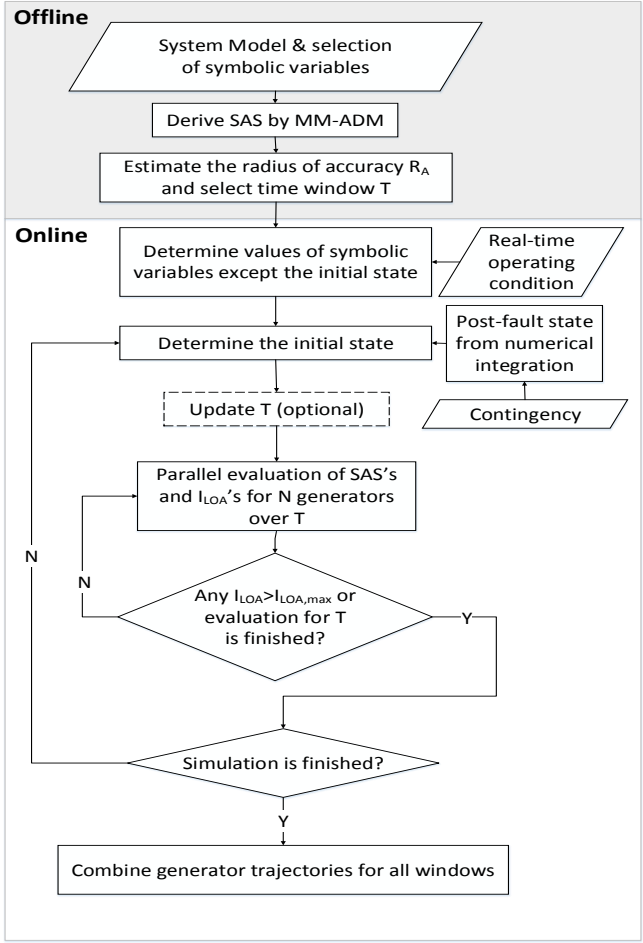

Fig. 8. Offline and online procedure diagram

This section summarizes the proposed approach using MM-ADM for power system simulation. It has an offline stage to derive the SAS and an online stage to evaluate the SAS as shown in Fig. 8.

### A. Offline Stage

A SAS is derived for each generator with symbolic variables from, e.g., one of these two groups:

*Group-1*: Time, the initial state, and the system operating condition (e.g. generator outputs and loads)

*Group-2*: *Group-1* plus selected symbolized elements in the system admittance matrix or the reduced admittance matrix having only generator internal nodes.

*Group-1* assumes a specific post-contingency system topology (i.e. the admittance matrix) but relaxes the system operating condition such as to enable one SAS to simulate for multiple loading conditions. *Group-2* additionally relaxes selected elements in the admittance matrix and hence enables one SAS to simulate for multiple contingencies. Other symbolic variables can also be added for any undetermined parameters, but the more symbolic variables, the more complex expression of the SAS. All SAS's derived in the offline stage will be saved in storage for online use.

The offline stage also needs to estimate $R_A$ using (35) from the system model and accordingly choose an appropriate length for the repeating time window, i.e. $T<R_A$.

### B. Online Stage

For a specific contingency, this stage evaluates the corresponding SAS of every generator consecutively for each $T$ until making up the expected simulation period. The first time window needs to know the initial state of the system, which can be obtained from numerical integration only for the fault-on period until the fault is cleared. Starting from the second time window, the initial state takes the last state of the previous time window. Because of the independence between the SAS's of different generators, evaluations of their SAS's can be performed in parallel on concurrent computers.

An adaptive time window $T$ may be applied. To do so, a $T$ smaller than the $R_A$ estimated for typical post-fault states is chosen as the initial time window. The aforementioned indicator $I_{LOA}$ will be calculated for each time window. If it is found hitting a preset limit at any time $t$ before the current window $T$ ends, the SAS evaluation for the following window of $T$ will start immediately. Also, it is optionally to re-evaluate $R_A$ using the initial state taken at $t$ to update the length of $T$ if it exceeds $R_A$. Thus, the actual time windows could be adaptively changed to ensure accuracy of the SAS.

## V. CASE STUDY ON THE IEEE 39-BUS SYSTEM

IEEE 10-generator, 39-bus system (as shown in Fig. 9) is used to validate the two-stage SAS-based approach for power system simulation. The original inertias are listed in Table III. The generator 39 has a much larger inertia than the others. The SAS with $N=3$ is derived for each generator. Then from (35), $R_A$ for each generator is estimated. Except for generator 39, all generators have a small $R_A$ around 0.04s, which is expected from the analysis on Fig. 4. With $T=0.04s$, rotor angles relative to generator 39 following a three-phase fault at bus 2 cleared by tripping the line 2-25 after 4 cycles are given from both RK-4 and the SAS as shown in Fig. 10 and Fig. 11, which match well.

According to (35), if $T$ is expected to be 0.2s and $I_{LOA,max}$ is selected as 3 rad/s, the minimum inertias for all generators are estimated as listed in the 3rd column of Table III. For instance, Table IV gives parameters for estimating the $H_{min}$ for generator 30 by selecting bus 39 as the reference node. Pick the largest $H_{min}$ value, i.e. 106, as the minimum expected inertia for all generators. Some random variations are added to create $H_{min}$ values for all generators, e.g. those in the last column of Table III. Consider the same contingency as above. Fig. 12 and Fig. 13 respectively give 4-second post-fault simulation results from the RK4 and the SAS on rotor angles relative to generator 39.

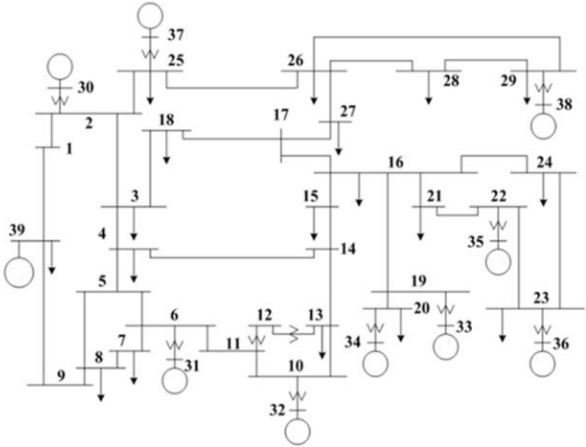

Fig. 9. IEEE 10-generator 39-bus system

TABLE III
INERTIAS AND LIMITS OF 10 GENERATORS (UNIT: S)

| Bus No. | Original H | $H_{min}$ for T=0.2s | $H_{min}$ used for T=0.2s |
|---|---|---|---|
| 30 | 4.2 | 106 | 106 |
| 31 | 3.03 | 54 | 109 |
| 32 | 3.58 | 47 | 105 |
| 33 | 2.86 | 21 | 110 |
| 34 | 2.6 | 10 | 113 |
| 35 | 3.48 | 25 | 104 |
| 36 | 2.64 | 14 | 107 |
| 37 | 2.43 | 53 | 111 |
| 38 | 3.45 | 18 | 110 |
| 39 | 50 | 50 | 114 |

TABLE IV
PARAMETERS (PU) OF BUS 30 AGAINST BUS 39

| $\theta_{12}$ | $Y_{12}$ | $G_1$ | $G_2$ | $\delta_1(0)$ | $\dot{\delta}_1(0)$ |
|---|---|---|---|---|---|
| 1.5458 | 7.4753 | 2.2361 | 1.5907 | -0.565 | -10.908 |
| $\delta_2(0)$ | $\dot{\delta}_2(0)$ | $P_{m1}$ | $P_{m2}$ | $E_1$ | $E_2$ |
| 0.0563 | 1.5480 | 2.5000 | 10.0573 | 1.0566 | 1 |

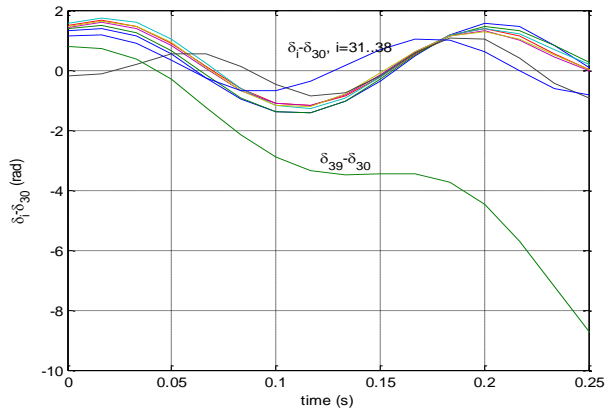

Fig. 10. Rotor angles from RK-4

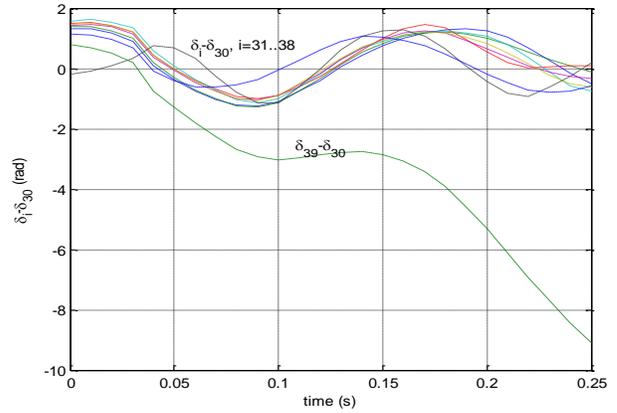

Fig. 11. Rotor angles from the SAS with T=0.04s

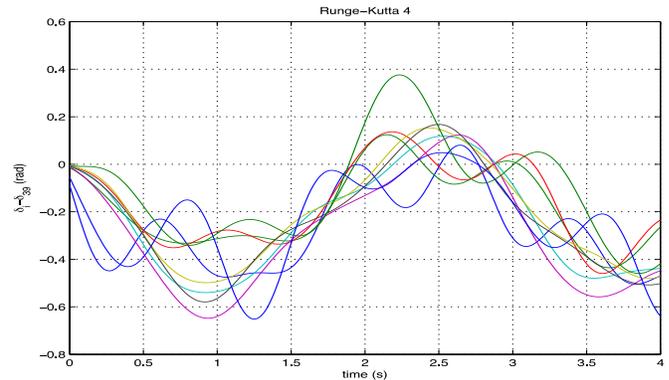

Fig. 12. IEEE 10-machine 39-bus system post-fault numerical integration

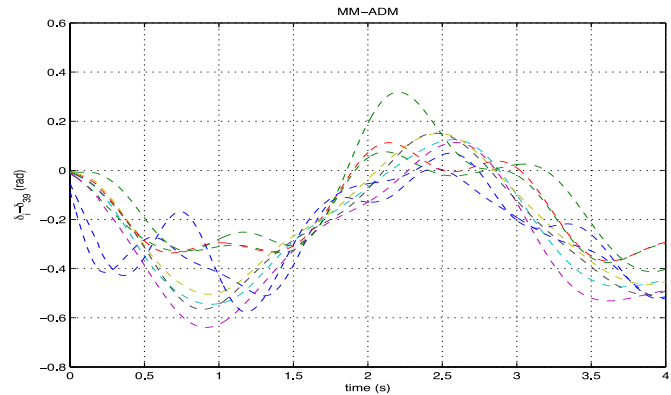

Fig. 13. IEEE 10-machine 39-bus system post-fault MM-ADM results

Because those $H_{min}$'s listed in table III are derived from one specific initial state, it is necessary to validate whether those $H_{min}$'s necessarily lead to a $R_A>0.2s$ for each generator thru the entire simulation period. According to Fig. 12 and Fig. 13, the worst initial state happens when $|\dot{\delta}(0)|$ reaches its maximum. Utilizing the data provided by Fig. 12, $|\dot{\delta}_{i-39}|_{max}$ for each machine is listed in Table V, among which the largest is 2.1rad/s about generator 30 happening at t=3.8667s. Apply the data at t=3.8667 as the initial state to (35), $R_A$ is solved as 0.35s>0.2s. Therefore the $H_{min}$'s in Table III are legitimate over the whole simulation period.

To better illustrate the accuracy of the SAS from MM-ADM, the comparison between RK4 and SAS results of 3 selected generators are shown in Fig. 14, Fig. 15 and Fig. 16, which match well.

TABLE V
$|\dot{\delta}_{i-39}|_{\text{MAX}}$ (RAD/S) WITHIN A 4S SIMULATION PERIOD

| $i$ | 30 | 31 | 32 | 33 | 34 | 35 | 36 | 37 | 38 |
|---|---|---|---|---|---|---|---|---|---|
| $|\dot{\delta}_{i-39}|_{\max}$ | 2.1 | 1.6 | 1.7 | 1.7 | 1.8 | 1.6 | 1.6 | 1.9 | 1.4 |

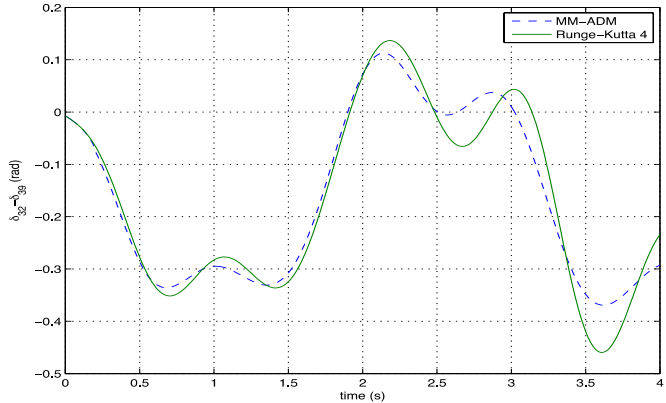

Fig. 14. $\delta_{32} - \delta_{39}$ Comparison between MM-ADM and Runge-Kutta 4

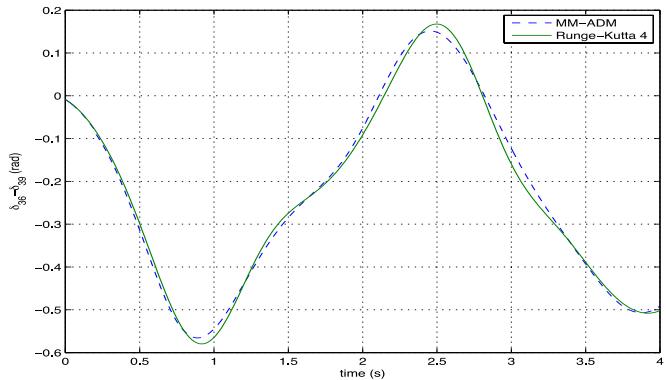

Fig. 15. $\delta_{36} - \delta_{39}$ Comparison between MM-ADM and Runge-Kutta 4

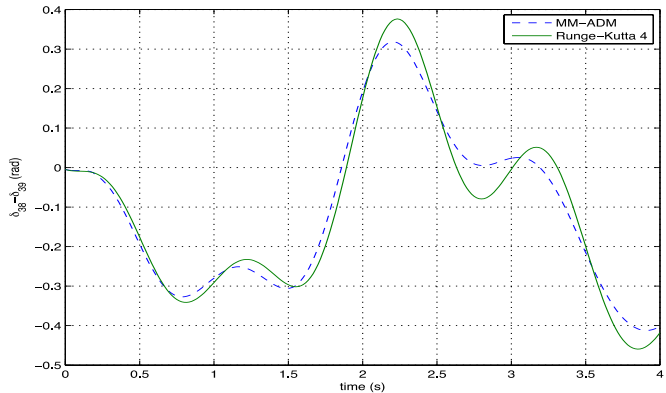

Fig. 16. $\delta_{38} - \delta_{39}$ Comparison between MM-ADM and Runge-Kutta 4

To demonstrate the time performance the proposed SAS-based approach, the following three cases are tested in the environment of MATLAB and MAPLE:

*Case-A:* only symbolizing time $t$, $\Delta\delta_i(0)$ and $\Delta\dot{\delta}_i(0)$, i.e. for one specific simulation.

*Case-B:* also symbolizing the reduced admittance matrix about 10 generator EMFs, i.e. for simulating different faults under one specific loading condition. The magnitude and angle of each element in reduced admittance matrix will be symbolized separately. Therefore, there will be two symmetric symbolic 10×10 matrices for magnitudes and angles.

*Case-C:* adding 10 symbolic variables for generators' mechanical powers to *Case-B* to make the SAS be also good for simulating different faults and loading conditions with each load modeled as a constant impedance.

The offline stage is implemented in MAPLE and the online stage is performed in MATLAB. Table VI lists the test results on time spent in the online and offline stages. For the online stage, we assume 3 points of the SAS to be evaluated for each window $T$: one at the middle point of $T$ and the other two near the end of $T$ to estimate the initial state for the next $T$. Thus, evaluations of SAS are given at about 0.1s intervals. That minimizes the computation with online SAS evaluation but still provides adequate information for judging transient stability.

TABLE VI
TIME PERFORMANCE

|  | Case-A | Case-B | Case-C |
|---|---|---|---|
| Offline time (s) | 218.11 | 392.40 | 395.64 |
| Online time (s) | 0.0045 | 0.0140 | 0.0140 |
| Ratio of $T$ to online time | 44 | 14 | 14 |

From the comparison between *Case-A* and *Case-B*, symbolizing the admittance matrix makes the offline and online computations take 80% and 210% longer, respectively, because 110 additional symbolic variables for the admittance matrix are added for *Case-B*. However, after adding more symbolic variables on the system operating condition, *Case-C* takes almost the same offline and online times as *Case-B*. The reason is that the symbolic variables on generations and loads are only involved in operations with addition, neither multiplication nor any sine function. For these three cases, the online computation times are 1/44, 1/14 and 1/14 of real time.

## VI. CONCLUSION

This paper has proposed a new approach for power system time-domain simulation, which is based on the SAS of power system DAEs derived by a MM-ADM based method. A two-stage scheme, i.e. offline SAS solving and online SAS evaluation, was presented to minimize the online computational burden.

## VII. REFERENCES

[1] G. Adomian, Solving Frontier Problems of Physics: The Decomposition Method, Kluwer Academic, Dordrecht, 1994.
[2] J. Biazar, S.M. Shafiof "A simple algorithm for calculating Adomian polynomials," *Int. J. Contemp. Math. Sciences*, vol. 2, no. 20, pp. 975–982, 2007.
[3] H. Jafari, C.M. Khalique & M.Nazari, "Application of the Laplace decomposition method for solving linear and nonlinear fractional diffusion-wave equations," Applied Mathematics Letter, vol. 24, pp. 1799–1805, Jan. 2011.
[4] J.S. Duan, R. Rach & A.M. Wazwaz, "A reliable algorithm for positive solutions of nonlinear boundary value problems by the multistage Adomian decomposition method," Applied Mathematical Modelling, vol. 37, no. 20, pp. 8687–8708, Jan.

[5] A. Wazwaz, "The modified decomposition method and Pade approximants for solving the Thomas Fermi equation," Applied Mathematics and Computation, vol. 105, no. 1, pp. 11–19, Jun. 1999.

[6] P.M. Anderson and A.A. Fouad, *Power system control and stability*. 2nd ed., New York: Wiley Interscience, 2003, pp.13-66.

[7] N.I. Razali, M.S. Chowdhury & W.Asrar, "The Multistage Adomian Decomposition Method for Solving Chaotic Lu System," *Middle-east J. Scientific Research*, vol. 13, pp. 43–49, Jan. 2013.

[8] G. Adomian, "On the convergence region for decomposition solutions,"*J. Computational Applied Mathematics*, vol. 11, pp. 379–380, 1984.

[9] J.S. Duan, et al, "New higher-order numerical one-step methods based on the Adomian and the modified decomposition methods," *Applied Mathematics and Computation*, vol. 218, pp. 2810–2828, 2011.

[10] O.T. Kolebaje, O.L. Akinyemi & R.A. Adenodi, "On the application of the multistage laplace adomian decomposition method with pade approximation to the rabinovich-fabrikant system," *Advances in Applied Science Research*, vol. 4, no. 3, pp. 232–243, 2013.

[11] S. M. Holmquist, " An examination of the effectiveness of the Adomian decomposition method in fluid dynamic applications," Ph.D. dissertation, Dept. Math., Univ. of Central Florida, Orlando, FL, 2007.

[12] A.S. Arife, S.T. Korashe & A.Yildirim, "Laplace Adomian decomposition method for solving a model chronic Myelogenous Leukemia CML and T Cell interaction," *World Applied Sciences J.*, vol. 13, no. 4, pp. 756–761, 2011.

[13] J.S. Duan, R. Rach, D.Baleanu & A. Wazwaz, "A review of the Adomian decomposition method and its applications to fractional differential equations," *Commun. Frac. Calc.*, vol. 3, no. 2, pp. 73–99, 2012.

[14] A. Wazwaz, "A reliable modifcation of Adomian decomposition method,"*Applied Mathematics and Computation*, vol. 102, pp. 77–86, 1999.

[15] S.A. Khuri, "A Laplace decomposition algorithm applied to a class of nonlinear differential equations," *J. Applied Mathematics*, vol. 1, pp. 141–155, Jun. 2001.

[16] L.A. Bougoffa, R.C. Rach & S.El-manouni, "A convergence analysis of the Adomian decomposition method for an abstract Cauchy problem of a system of firstorder nonlinear differential equations," *Int. J. Comput. Mathematics*, vol. 2, no. 90, pp. 360–375, Aug. 2012.